\documentclass{article}
\usepackage[english]{babel}

\usepackage{theorem}
\usepackage{latexsym}
\usepackage{amsfonts}
\usepackage{amssymb}
\usepackage{multicol}

\newenvironment{proof}{\noindent{\it Proof. }}{\fbox{}\\}
\newtheorem{theorem}{Theorem}
\newtheorem{proposition}{Proposition}
\newtheorem{lemma}{Lemma}
\newtheorem{corollary}{Corollary}
{\theorembodyfont {\rmfamily}

\newtheorem{definition}{Definition}
\newtheorem{example}{Example}}

\begin{document}

\title{\bf{Free-differentiability conditions on the free-energy
function implying large deviations}}

\author{Henri Comman\thanks{Department of Mathematics, University
of Santiago de Chile, Bernardo O'Higgins 3363, Santiago, Chile.
E-mail: hcomman@usach.cl}}

\date{}

\maketitle

\abstract{
 Let $(\mu_{\alpha})$ be a net of Radon sub-probability measures on $\mathbb{R}$, and
  $(t_{\alpha})$ be a net in $]0,+\infty[$ converging to $0$. Assuming that the generalized
  log-moment generating function $L(\lambda)$ exists for all
$\lambda$ in a nonempty open interval $G$, we give conditions on
the left or right derivatives of $L_{\mid G}$, implying   vague
(and thus narrow when $0\in G$) large deviations. The rate
function (which can be nonconvex)  is obtained as
     an abstract Legendre-Fenchel transform. This allows us to strengthen the
      G\"{a}rtner-Ellis theorem by removing the usual differentiability assumption.
      A related question of R. S. Ellis is solved.}

\section{Introduction}

 Let $(\mu_{\alpha})$ be a net of
 Radon sub-probability measures on a Hausdorff topological space
 $X$, and
  $(t_{\alpha})$ be a net in $]0,+\infty[$ converging to $0$. Let
 $\mathcal{B}(X)$ (resp. $\mathcal{C}(X)$) denote the set of
$[-\infty,+\infty[$-valued Borel measurable (resp. continuous)
  functions on $X$.
 For each $h\in\mathcal{B}(X)$,
 we define \[\underline{\Lambda} (h)=\log\liminf
\mu_{\alpha}^{t_{\alpha}}(e^{h/t_{\alpha}})\] and
\[\overline{\Lambda} (h)=\log\limsup
\mu_{\alpha}^{t_{\alpha}}(e^{h/t_{\alpha}})\] where
$\mu_{\alpha}^{t_{\alpha}}(e^{h/t_{\alpha}})$ stands for $(\int_X
e^{h(x)/t_{\alpha}}\mu_{\alpha}(dx))^{t_{\alpha}}$, and   write
$\Lambda(h)$ when both expressions are equal.
 When
$X=\mathbb{R}$, for each pair of reals $(\lambda,\nu)$, let
$h_{\lambda,\nu}$ be  the function   defined on $X$ by
$h_{\lambda,\nu}(x)=\lambda x$ if $x\le 0$ and
$h_{\lambda,\nu}(x)=\nu x$ if $x\ge 0$ (we write simply
$h_\lambda$ in place of $h_{\lambda,\lambda}$).  For each real
$\lambda$,  we put
 $L(\lambda)=\Lambda(h_\lambda)$
when $\Lambda(h_\lambda)$ exists.

A well-known problem of large deviations in $\mathbb{R}$ (usually
stated for sequences of probability measures) is the following:
assuming that $L(\lambda)$ exists and is finite for all $\lambda$
in an open interval $G$ containing $0$, and that the map $L_{\mid
G}$  is not differentiable on $G$, what conditions on $L_{\mid G}$
do imply large deviations, and with which rate function?

In relation with this problem, R. S. Ellis posed the following
question (\cite{ell1}): assuming that $\Lambda(h_{\lambda,\nu})$
exists and is finite for all $(\lambda,\nu)\in\mathbb{R}^2$,  what
conditions on the functional $\Lambda_{\mid \{
h_{\lambda,\nu}:(\lambda,\nu)\in\mathbb{R}^2\}}$ do
 imply large deviations with rate function
$J(x)=\sup_{(\lambda,\nu)\in\mathbb{R}^2}\{h_{\lambda,\nu}(x)-\Lambda(h_{\lambda,\nu})\}$
for all $x\in X$ ?

In this paper, we solve the above  problem by giving conditions on
$L_{\mid G}$ involving only its left and right derivatives; the
rate function is obtained as an abstract Legendre-Fenchel
transform ${\Lambda_{\mid \mathcal{S}}}^*$, where $\mathcal{S}$
can be any set in $\mathcal{C}(X)$ containing
$\{h_\lambda:\lambda\in G \}$ (Theorem \ref{open-problem}). When
$\mathcal{S}=\{h_\lambda:\lambda\in G\}$, we get a strengthening
of the
 G\"{a}rtner-Ellis theorem by removing the usual differentiability
 assumption (Corollary \ref{GE}). The answer to the Ellis question is obtained with
$\mathcal{S}=\{h_{\lambda,\nu}:(\lambda,\nu)\in\mathbb{R}^2\}$
 (Corollary \ref{open-question-Ellis}).

 The techniques
used  are refinements of those developed in previous author's
works (\cite{com1}, \cite{com2}), where variational forms for
$\underline{\Lambda}(h)$ and $\overline{\Lambda}(h)$ are obtained
with $h\in\mathcal{B}(X)$  satisfying the usual Varadhan's tail
condition ($X$ a general space). We consider here the set
$\mathcal{C}_{\mathcal{K}}(X)$ of elements  $h$ in
$\mathcal{C}(X)$ for which $\{y\in X:e^{h(x)}
 -\varepsilon\le e^{h(y)}\le e^{h(x)} +\varepsilon\}$ is compact for all $x\in
 X$ and $\varepsilon>0$ with $e^{h(x)}>\varepsilon$.
 The first step is Theorem
\ref{compact-case}, which establishes
  that for any
  $\mathcal{T}\subset\mathcal{C}_{\mathcal{K}}(X)$, and under suitable conditions
  (weaker than vague large deviations), there exist some reals
  $m,M$ such that
\[
\Lambda(h)=\sup_{x\in \{m\le h\le M\}}\{h(x)-l_1(x)\}\ \ \ \ \ \ \
\textnormal{for all}\ h\in\mathcal{T},
\]
where $l_1(x)= -\log\inf\{\liminf\mu_{\alpha}^{t_{\alpha}}(G):x\in
G\subset X, G\textnormal{\ open}\}$ for all $x\in X$; in
particular, $\Lambda(h)$ exists and has the same form as when
large deviations hold. Note that when $X=\mathbb{R}$ and
$\mathcal{T}=\{h_\lambda:\lambda\in G\}$ with $0\not\in G$, then
the $\sup$ in the above expression can be taken on a compact set
(if $0\in G$, this follows from the exponential tightness). It
turns out  that any subnet of $(\mu_{\alpha}^{t_{\alpha}})$ has a
subnet $(\mu_{\gamma}^{t_{\gamma}})$ satisfying the above
conditions. The second step consists then in applying Theorem
\ref{compact-case} with $X=\mathbb{R}$,
$\mathcal{T}=\{h_\lambda:\lambda\in G\}$ and all these subnets.
 More precisely,  we show that if $x$ is the left or right
derivative of $L$ at some point $\lambda_x\in G$, then
$l^{(\mu_{\gamma}^{t_{\gamma}})}_1(x)\le\lambda_x x-L(\lambda_x)$,
  whence
\begin{equation}\label{intro-eq2}
l^{(\mu_{\gamma}^{t_{\gamma}})}_1(x)\le {L_{\mid G}}^*(x)
\end{equation}
(Proposition \ref{lem-x+}).
  Let $\mathcal{S}$ be any set in $\mathcal{C}(X)$ containing $\{h_\lambda:\lambda\in G\}$,
 and assume that $\Lambda(h)$ exists for all $h\in\mathcal{S}$.
 It is easy to see  that
\begin{equation}\label{intro-eq3}
{L_{\mid G}}^*\le{\Lambda_{\mid \mathcal{S}}}^*\le
l^{(\mu_{\gamma}^{t_{\gamma}})}_0\le
l^{(\mu_{\gamma}^{t_{\gamma}})}_1,
\end{equation}
where  $l^{(\mu_{\gamma}^{t_{\gamma}})}_0(x)=
-\log\inf\{\limsup\mu_{\gamma}^{t_{\gamma}}(G):x\in G\subset X, G\
\textnormal{open}\}$ for all $x\in X$. Putting together
(\ref{intro-eq2}) and (\ref{intro-eq3}) give
\begin{equation}\label{intro-eq4}
{L_{\mid G}}^*(x)={\Lambda_{\mid
\mathcal{S}}}^*(x)=l^{(\mu_{\gamma}^{t_{\gamma}})}_0(x)=l^{(\mu_{\gamma}^{t_{\gamma}})}_1(x)
\end{equation}for all $x$ in the image of the left (resp. right) derivative of $L_{\mid
G}$;  consequently, if the set of these images contains
$\{{\Lambda_{\mid \mathcal{S}}}^*<+\infty\}$, then
$(\mu_{\gamma}^{t_{\gamma}})$ satisfies a  vague (narrow if $0\in
G$) large deviation principle with powers $(t_{\gamma})$ and rate
function ${\Lambda_{\mid \mathcal{S}}}^*$, which moreover
coincides with ${L_{\mid G}}^*$ on its effective domain. By
compactness and Hausdorffness arguments, we conclude  that the
same result holds for the net $(\mu_{\alpha}^{t_{\alpha}})$.
Furthermore, $\{{\Lambda_{\mid \mathcal{S}}}^*<+\infty\}$ can be
replaced by its interior, when ${\Lambda_{\mid \mathcal{S}}}^*$ is
proper convex and lower semi-continuous, which is the case when
$\mathcal{S}=\{h_\lambda:\lambda\in G\}$; this allows us to
improve a strong version of G\"{a}rtner-Ellis theorem given by O'
Brien.

Various generalizations  are given in order to get large
deviations  with a rate function coinciding with ${\Lambda_{\mid
\mathcal{S}}}^*$ and ${L_{\mid G}}^*$ only on its effective
domain.
 Note that
all our results hold for general nets of sub-probability measures
and powers.

The paper is organized as follows. Section \ref{section-prelim}
fixes the  notations and recall some results on large deviations
and convexity; Section \ref{section-variational-forms} deals with
the variational forms of the functionals $\Lambda$; Section
\ref{section-real-case} treats the case $X=\mathbb{R}$.

\section{Preliminaries}\label{section-prelim}

  Throughout the paper, the notations $\underline{\Lambda}$,
  $\overline{\Lambda}$, $\Lambda$, $l_0$, $l_1$
refer to the net  $(\mu_{\alpha}^{t_\alpha})$.
 We shall write  $l_1^{(\mu_{\beta}^{t_\beta})}$ when  in the definition of
 $l_1$,
  $(\mu_{\alpha}^{t_\alpha})$ is replaced by the subnet
  $(\mu_{\beta}^{t_\beta})$.
  We do not make such distinction
 for the map $\Lambda$, since it does not depend on the subnet along
 which the limit is taken. We recall  that $l_0$
and $l_1$ are lower semi-continuous functions.

 \begin{definition}
 \ \
 \begin{itemize}
 \item[(a)]
$(\mu_{\alpha})$ satisfies a \textit{(narrow) large deviation
principle} with powers $(t_{\alpha})$  if there exists a
$[0,+\infty]$-valued lower semi-continuous function $J$  on $X$
such that
\begin{equation}\label{defi-upper-ldp}
\limsup\mu_{\alpha}^{t_{\alpha}}(F)\le\sup_{x\in F}e^{-J(x)}\ \ \
\ \ \ \ \textnormal{for all closed $F\subset X$}
\end{equation}and
\[
\sup_{x\in G}e^{-J(x)}\le
 \liminf\mu_{\alpha}^{t_{\alpha}}(G)\ \ \
\ \ \ \ \textnormal{for all open $G\subset X$};\] $J$ is a
\textit{rate function} for $(\mu_{\alpha}^{t_{\alpha}})$, which is
said to be tight when it has compact level sets. When "closed" is
replaced by "compact"  in (\ref{defi-upper-ldp}), we say that a
\textit{vague large deviation principle} holds.
\item[(b)] $(\mu_\alpha)$ is \textit{exponentially tight with respect to
$(t_\alpha)$}  if for each $\varepsilon>0$ there exists a compact
set $K_\varepsilon\subset X$  such that
$\limsup\mu_{\alpha}^{t_{\alpha}}(X\verb'\'K_\varepsilon)<\varepsilon$.
\end{itemize}
\end{definition}

The following results are well-known for a net
$(\mu_\varepsilon^{\varepsilon})_{\varepsilon>0}$, with
$\mu_\varepsilon$  a Radon  probability measure (\cite{dem}); it
is easy to see that the proofs work also for  general nets of
sub-probability measures  and powers.

\begin{lemma}\label{lem-prelim}
\ \
\begin{itemize}
\item[(a)] Let $X$ be locally compact Hausdorff. Then, $(\mu_{\alpha})$ satisfies a
vague large deviation principle with powers $(t_{\alpha})$ if and
only if  $l_0=l_1$. In this case, $l_0$ is the rate function.
\item[(b)]  If $(\mu_{\alpha})$ satisfies a vague large deviation principle with
powers $(t_{\alpha})$, and $(\mu_\alpha)$ is exponentially tight
with respect to $(t_\alpha)$, then $(\mu_{\alpha})$ satisfies a
large deviation principle with same powers and same rate function.
\end{itemize}
\end{lemma}

A capacity on $X$ is a map $c$ from the powerset of $X$ to
$[0,+\infty]$ such that:
\begin{itemize}
\item[(i)] $c(\emptyset)=0$.
\item[(ii)] $c(Y)=\sup\{c(K):K\subset Y,K\textnormal{\ compact}\}\ \ \ $ for all
$Y\subset X$.
\item[(iii)] $c(K)=\inf\{c(G):K\subset G\subset X, G\ \textnormal{open}\}\ \ \ \ $ for
all compact $K\subset X$.
\end{itemize}
 The
vague topology on the set of capacities is the coarsest topology
for which the maps $c\rightarrow c(Y)$ are upper (resp. lower)
semi-continuous for all compact (resp. open) $Y\subset X$.  Let
$\Gamma(X,[0,1])$ denote the set of $[0,1]$-valued capacities on
$X$ provided with the vague topology, and note that
$(\mu_{\alpha}^{t_{\alpha}})$ is a net in $\Gamma(X,[0,1])$. For
each $[0,+\infty]$-valued lower semi-continuous function $l$ on
$X$, we associate the element $c_l$ in $\Gamma(X,[0,1])$ defined
by $c_l(Y)=\sup_{x\in Y} e^{-l(x)}$ for all $Y\subset X$. We refer
to \cite{ob2}  for the first assertion  in the following lemma;
the second one is  the mere transcription of the definition of a
vague large deviation principle in terms of capacities.

\begin{lemma}\label{lem-capa}
\ \
\begin{itemize}
\item[(a)] If  $X$ is  locally compact Hausdorff, then
  $\Gamma(X,[0,1])$ is
a compact Hausdorff space.
\item[(b)] $(\mu_{\alpha})$ satisfies a vague large deviation principle
with powers $(t_{\alpha})$ and rate function $J$ if and only if
$(\mu_{\alpha}^{t_{\alpha}})$ converges to $c_{J}$ in
$\Gamma(X,[0,1])$.
\end{itemize}
\end{lemma}

For any $[-\infty,+\infty]$-valued (not necessary convex) function
$f$ defined on some topological  space, we put
$\mathcal{D}\textnormal{om}(f)=\{f<+\infty\}$ (the so-called
effective domain), and denote by
  $\textnormal{int}\mathcal{D}\textnormal{om}(f)$ (resp.
  $\textnormal{bd}\mathcal{D}\textnormal{om}(f)$) the   interior (resp.
  boundary) of $\mathcal{D}\textnormal{om}(f)$. The range of $f$ is denoted by
  $\textnormal{ran}f$.

 A $[-\infty,+\infty]$-valued convex function $f$ on $\mathbb{R}$ is said to
   be proper if $f$ is $]-\infty,+\infty]$-valued and takes a finite value on at least one point.
 The Legendre-Fenchel transform $f^*$ of $f$
  is defined by $f^*(x)=\sup_{\lambda\in\mathbb{R}}\{\lambda
  x-f(\lambda)\}$ for all  $x\in \mathbb{R}$; note that $f^*$ is  convex lower
semi-continuous, and  proper when $f$ is proper.
   Let $I\subset \mathbb{R}$ be a nonempty interval, and $f_{\mid I}$ be a
   $]-\infty,+\infty]$-valued
   convex function on $I$. We denote by $\widehat{f_{\mid I}}$ the convex
   function on $\mathbb{R}$ which coincides with $f_{\mid I}$ on $I$, and takes the value
   $+\infty$ out $I$;
    in this case we write simply
 ${f_{\mid I}}^*$ in place of ${\widehat{f_{\mid I}}}^*$.
    The left and  right derivatives of $f_{\mid I}$ at some point
    $x\in\mathcal{D}\textnormal{om}(f_{\mid I})$
     are denoted by ${f_{\mid I}}'_-(x)$ and
  ${f_{\mid I}}'_+(x)$ respectively.
A proper convex function $f$ on $\mathbb{R}$ is said to be
essentially smooth if
$\textnormal{int}\mathcal{D}\textnormal{om}(f)\neq\emptyset$,
  $f$ is differentiable on
 $\textnormal{int}\mathcal{D}\textnormal{om}(f)$,  and
 $\lim|f'(x_n)|=+\infty$ for all sequences $(x_n)$ in
 $\textnormal{int}\mathcal{D}\textnormal{om}(f)$ converging to
 some $x\in\textnormal{bd}\mathcal{D}\textnormal{om}(f)$
 (\cite{roc}).

If
 $L(\lambda)$ exists and is finite for all  $\lambda$ in a nonempty open interval
 $G$, then $L_{\mid G}$ is convex; if moreover
$0\in G$, then $(\mu_\alpha)$ is exponentially tight with respect
to $(t_\alpha)$. If $L(\lambda)$ exists for all reals $\lambda$,
then
  $L$ is a
$[-\infty,+\infty]$-valued convex function on $\mathbb{R}$; if
moreover $0\in\textnormal{int}\mathcal{D}\textnormal{om}(L)$, then
$L$ is proper
 (the proof of these facts is obtained by modifying suitably the one of
 Lemma 2.3.9 in
    \cite{dem}).

\begin{lemma}\label{conv}
Let $f$ be a proper convex lower semi-continuous function on
$\mathbb{R}$.
 Then,
\[\inf_{y\in G}f(y)=\inf_{y\in
G\cap\textnormal{int}\mathcal{D}\textnormal{om}(f)}f(y)\] for all
open sets $G\subset\mathbb{R}$.
\end{lemma}

\begin{proof}
 Let $G$ be an open subset of $\mathbb{R}$. If
 $G\cap\mathcal{D}\textnormal{om}(f)=\emptyset$,
 then the conclusion holds trivially ($\inf\emptyset=+\infty$ by convention). Assume that
 $G\cap\mathcal{D}\textnormal{om}(f)\neq\emptyset$. By Corollary
 6.3.2 of \cite{roc},
 $G\cap\textnormal{int}\mathcal{D}\textnormal{om}(f)\neq\emptyset$.
 By Theorem VI.3.2 of \cite{ell2}, for each $x\in\mathcal{D}\textnormal{om}(f)$
  we can find a sequence $(x_n)$ in
$\textnormal{int}\mathcal{D}\textnormal{om}(f)$ converging to $x$
and such that $\lim f(x_n)=f(x)$, which implies
$\inf_{G\cap\mathcal{D}\textnormal{om}(f)}f=\inf_{G\cap\textnormal{int}
\mathcal{D}\textnormal{om}(f)}f$, and  the lemma is proved since
$\inf_{G\cap\mathcal{D}\textnormal{om}(f)}f=\inf_G f$.
\end{proof}

\section{Variational forms for $\Lambda$ on $\mathcal{C}_{\mathcal{K}}(X)$}
\label{section-variational-forms}

We begin by defining a notion, which will appear as a key
condition in the sequel; it is nothing else but a uniform version
of  the tail condition in Varadhan's theorem.

\begin{definition}
We say that a set $\mathcal{T}\subset\mathcal{B}(X)$
\textit{satisfies the tail condition} for
$(\mu_{\alpha}^{t_{\alpha}})$ if for each $\varepsilon>0$, there
exists a real $M$ such that
\[\limsup\mu_{\alpha}^{t_{\alpha}}(e^{h/t_{\alpha}}1_{\{h>M\}})<\varepsilon
\ \ \ \ \ \ \textnormal{for all $h\in\mathcal{T}$}.
\]
\end{definition}

For each  $h\in\mathcal{B}(X)$, each $x\in X$ and each
$\varepsilon>0$, we put $F_{e^{h(x)},\varepsilon}=\{y\in
X:e^{h(x)}
 -\varepsilon\le e^{h(y)}\le e^{h(x)} +\varepsilon\}$ and
 $G_{e^{h(x)},\varepsilon}=\{y\in X:e^{h(x)}
 -\varepsilon<e^{h(y)}<e^{h(x)} +\varepsilon\}$.
The following expressions are known when $(\mu_\alpha)$ is a net
of probability measures, and when $\mathcal{T}$ has only one
element, say $h$ (see \cite{com1} and \cite{com2}  for the first
and the second assertion, respectively). The proofs reveal that
the constant $M$ comes from the above tail condition (assumed to
be satisfied by $h$),  so that the uniform versions for a general
$\mathcal{T}$ follow immediately; they moreover work as well for
the sub-probability case.

\begin{theorem}\label{lim}
Let $\mathcal{T}\subset\mathcal{B}(X)$ satisfying the tail
condition for $(\mu_{\alpha}^{t_{\alpha}})$. There is a real $M$
such that for each $h\in\mathcal{T}$,
\[
e^{\underline{\Lambda}(h)}=\liminf\sup_{x\in
 X,\varepsilon>0}\{(e^{h(x)}-\varepsilon)\mu_{\alpha}^{t_{\alpha}}(
G_{e^{h(x)},\varepsilon})\}=\lim_{\varepsilon\rightarrow
0}\liminf\sup_{x\in\{h\le M\}}\{e^{h(x)}
\mu_{\alpha}^{t_{\alpha}}(G_{e^{h(x)},\varepsilon})\}
\]
and
\[
e^{\overline{\Lambda}(h)}=\sup_{x\in X
,\varepsilon>0}\{(e^{h(x)}-\varepsilon)\limsup\mu_{\alpha}^{t_{\alpha}}(
G_{e^{h(x)},\varepsilon})\}=\sup_{x\in\{h\le
M\},\varepsilon>0}\{(e^{h(x)}-\varepsilon)\limsup\mu_{\alpha}^{t_{\alpha}}(
G_{e^{h(x)},\varepsilon})\}.
\]
In the above expressions, $G_{e^{h(x)},\varepsilon}$ can be
replaced by $F_{e^{h(x)},\varepsilon}$.
\end{theorem}

Part (a) of the following theorem shows that under conditions
strictly weaker than large deviations,  $\Lambda(h)$ exists and
has the same form as when large deviations hold, since in this
case the rate function coincides with $l_1$ (Lemma
\ref{lem-prelim}); it can be seen as a vague version of Varadhan's
theorem. Note that the hypothesis
$h\in\mathcal{C}_{\mathcal{K}}(X)$ cannot be dropped:
 consider a vague large deviation
principle for a net of probability measures with rate function
$J\equiv+\infty$,  take $h\equiv 0$ and get $\Lambda(h)=0$ and
$\sup_{X}\{h(x)-J(x)\}=-\infty$. Note also that the condition
$(ii)$ holds in particular when $(\mu_{\alpha}^{t_{\alpha}})$
converges in $\Gamma(X,[0,1])$.

\begin{theorem}\label{compact-case}
Let $\mathcal{T}\subset\mathcal{C}(X)$ with $X$ locally compact
Hausdorff, and assume that the following hold:
\begin{itemize}
\item[(i)] $\mathcal{T}$ satisfies the tail
condition for $(\mu_{\alpha}^{t_{\alpha}})$.
\item[(ii)]
$\limsup\mu_{\alpha}^{t_{\alpha}}(K)\le\liminf\mu_{\alpha}^{t_{\alpha}}(G)$
for each compact $K\subset X$ and each open $G\subset X$ with
$K\subset G$.
\item[(iii)] $\inf_{h\in\mathcal{T}}\overline{\Lambda}(h)>m$ for some real $m$.
\end{itemize}
 The following conclusions hold.
\begin{itemize}
\item[(a)] If  $\mathcal{T}\subset\mathcal{C}_{\mathcal{K}}(X)$, then
$\Lambda(h)$ exists for all $h\in\mathcal{T}$, and there is a real
$M$ such that
\begin{equation}\label{compact-case-eq1}
\Lambda(h)=\sup_{x\in\{m\le h\le M\}} \{h(x)-l_1(x)\}=\sup_{x\in
X} \{h(x)-l_1(x)\}\ \ \ \ \ \textnormal{for all
$h\in\mathcal{T}$}.
\end{equation}
\item[(b)] If  $(\mu_\alpha)$
is exponentially tight with respect to $(t_\alpha)$, then
$\Lambda(h)$ exists for all $h\in\mathcal{T}$, and there is a real
$M$ and a compact $K\subset X$ such that
\begin{equation}\label{compact-case-eq2}
\Lambda(h)=\sup_{x\in K\cap\{m\le h\le M\}}
\{h(x)-l_1(x)\}=\sup_{x\in X} \{h(x)-l_1(x)\}\ \ \ \ \
\textnormal{for all $h\in\mathcal{T}$}.
\end{equation}
\end{itemize}
\end{theorem}

\begin{proof}
Assume $\mathcal{T}\subset\mathcal{C}_{\mathcal{K}}(X)$.
 By $(i)$ and
 Theorem \ref{lim}, there is a real $M'$ such that for each
 $h\in\mathcal{T}$,
\begin{equation}\label{compact-case-eq3}
\sup_{x\in \{h\le M'+\log 2\}}e^{h(x)}e^{-l_1(x)}\le\sup_{x\in X
}e^{h(x)}e^{-l_1(x)}\le e^{\underline{\Lambda}(h)}\end{equation}
\[\le
e^{\overline{\Lambda}(h)}=\sup_{x\in\{h\le
M'\},\varepsilon>0}\{(e^{h(x)}-\varepsilon)\limsup\mu_{\alpha}^{t_{\alpha}}
(F_{e^{h(x)},\varepsilon})\}.
\]
Put $M=\log 2+M'$, and  suppose that
\[
\sup_{x\in \{h\le M\}}e^{h(x)}e^{-l_1(x)}+\nu<\sup_{x\in\{h\le
M'\},\varepsilon>0}\{(e^{h(x)}-\varepsilon)\limsup\mu_{\alpha}^{t_{\alpha}}
(F_{e^{h(x)},\varepsilon})\}
\] for some $h\in\mathcal{T}$ and some $\nu>0$. Then there exists $x_0\in\{h\le M'\}$
and  $\varepsilon_0>0$ with $e^{h(x_0)}>\varepsilon_0$ such that
\begin{equation}\label{compact-case-eq4}
\sup_{x\in \{h\le M\}}e^{h(x)}e^{-l_1(x)}<
(e^{h(x_0)}-\varepsilon_0-\nu)\limsup\mu_{\alpha}^{t_{\alpha}}(F_{e^{h(x_0)},\varepsilon_0}).
\end{equation}
By continuity and local compactness, for each $x\in
F_{e^{h(x_0)},\varepsilon_0}$, there exist some open sets $V_x$
and $V'_x$  satisfying $x\in V_x\subset\overline{V_x}\subset V'_x$
with $\overline{V_x}$ compact, and  such that
$e^{h(y)}>e^{h(x_0)}-\varepsilon_0-\nu$ for all $y\in V'_x$. Note
that
 $h(x)\le M$ for each  $x\in F_{e^{h(x_0)},\varepsilon_0}$,
since $e^{h(x_0)}+\varepsilon_0<2e^{M'}$. By
(\ref{compact-case-eq4}), for each $x\in
F_{e^{h(x_0)},\varepsilon_0}$, there exist some open sets $W_x$
and $W'_x$  satisfying $x\in W_x\subset\overline{W_x}\subset W'_x$
with $\overline{W_x}$ compact, and such that
\begin{equation}\label{compact-case-eq6}
e^{h(x)}\liminf\mu_{\alpha}^{t_{\alpha}}(W'_x)<(e^{h(x_0)}-\varepsilon_0-\nu)
\limsup\mu_{\alpha}^{t_{\alpha}}(F_{e^{h(x_0)},\varepsilon_0}).
\end{equation}
Put $G_x=W_x\cap V_x$ for all $x\in F_{e^{h(x_0)},\varepsilon_0}$.
 Since $F_{e^{h(x_0)},\varepsilon_0}$ is compact, there is a
finite set $A\subset F_{e^{h(x_0)},\varepsilon_0}$ such that
$F_{e^{h(x_0)},\varepsilon_0}\subset\bigcup_{x\in A} G_x$; thus,
for  some $x\in A$ we have
\[
(e^{h(x_0)}-\varepsilon_0-\nu)\limsup\mu_{\alpha}^{t_{\alpha}}(F_{e^{h(x_0)},\varepsilon_0})\le
e^{h(x)}\limsup\mu_{\alpha}^{t_{\alpha}}(G_x)\]
\[\le e^{h(x)}
\limsup\mu_{\alpha}^{t_{\alpha}}(\overline{W_x})\le e^{h(x)}
\liminf\mu_{\alpha}^{t_{\alpha}}(W'_x)
\]
(where the third inequality follows from $(ii)$),  which
contradicts
 (\ref{compact-case-eq6}). Therefore, all inequalities in
 (\ref{compact-case-eq3}) are equalities, that is  for each $h\in\mathcal{T}$,
 $\Lambda(h)$ exists and
 \[\Lambda(h)=\sup_{x\in \{h\le M\}}\{h(x)-l_1(x)\}=
 \sup_{x\in X}\{h(x)-l_1(x)\}= \sup_{x\in \{m\le h\le M\}}\{h(x)-l_1(x)\},\]
 (where the third equality follows from
 $(iii)$),
which proves $(a)$. For $(b)$, the above proof works verbatim
replacing $\{h\le M\}$ and $F_{e^{h(x_0)},\varepsilon_0}$ by
$\{h\le M\}\cap K$ and $F_{e^{h(x_0)},\varepsilon_0}\cap K$
respectively, where $K$ is some compact set  given by the
exponential tightness.
\end{proof}

The following definition extends the usual notion of
Legendre-Fenchel transform (when $X$ is a real topological vector
space and $\mathcal{S}$ its topological dual) and its
generalization proposed in \cite{ell1} (with $X=\mathbb{R}$ and
$\mathcal{S}=\{h_{\lambda,\nu}:(\lambda,\nu)\in\mathbb{R}^2\}$);
it coincides with our preceding notations since for
$\mathcal{S}=\{h_\lambda:\lambda\in G\}$ with $G$  a nonempty open
interval, we have
\[{L_{\mid G}}^*(x)=\sup_{\lambda\in\mathbb{R}}\{\lambda
x-\widehat{L_{\mid G}}(\lambda)\}=\sup_{\lambda\in G}\{\lambda
x-L(\lambda)\}=\sup_{\{h_\lambda:\lambda\in G\}}\{h_\lambda(
x)-\Lambda(h_\lambda)\}={\Lambda_{\mid S}}^*(x).\]
 In \cite{com1} (Corollary 2), we proved that for $X$ completely regular (not necessary Hausdorff),
   a rate function has always the form ${\Lambda_{\mid \mathcal{S}}}^*$,
   where $\mathcal{S}$ is
    any set in $\mathcal{C}(X)$ stable by translation,
   separating  suitably points and closed sets, and such that
   each $h\in\mathcal{S}$ satisfies the tail condition for
   $(\mu_{\alpha}^{t_{\alpha}})$;
    this is
    proved
   in \cite{com2} for $X$ normal Hausdorff and
   $\mathcal{S}$ the set of all bounded continuous
   functions on $X$ (this case  was known under exponential tightness hypothesis
   as a part of the conclusion
   of Bryc's theorem). We will identify in the next section others
   sets $\mathcal{S}$ for which the rate function is given by ${\Lambda_{\mid
   \mathcal{S}}}^*$.

\begin{definition}
Let $\mathcal{S}\subset\mathcal{B}(X)$ such that $\Lambda(h)$
exists for all $h\in \mathcal{S}$. The map ${\Lambda_{\mid
\mathcal{S}}}^*$ defined by
\[{\Lambda_{\mid\mathcal{S}}}^*(x)=\sup_{h\in\mathcal{S}}\{h(x)-\Lambda(h)\}
\ \ \ \ \ \ \ \textnormal{for all $x\in X$},\] is the
\textsl{abstract Legendre-Fenchel transform  of $\Lambda_{\mid
\mathcal{S}}$}.
\end{definition}

\section{The case  $X=\mathbb{R}$}\label{section-real-case}

In this section, we take $X=\mathbb{R}$ and   apply Theorem
\ref{compact-case}
 with
$\mathcal{T}=\{h_\lambda:\lambda\in G\}$ where $G$ is a nonempty
open interval.
 This allows us  to  compare  the values of $l^{(\mu_{\gamma}^{t_\gamma})}_1$ and  those of
${L_{\mid G}}^*$ on $\textnormal{ran}{L_{\mid
G}}_-'\cup\textnormal{ran}{L_{\mid G}}_+'$, where
$(\mu_{\gamma}^{t_\gamma})$ is a suitable  subnet of
$(\mu_{\alpha}^{t_\alpha})$  (Proposition \ref{lem-x+}). By means
of  a compactness argument, we then  derive sufficient conditions
for large deviations, involving only the left and right
derivatives of $L_{\mid G}$; the rate function is given by an
abstract  Legendre-Fenchel transform
${\Lambda_{\mid\mathcal{S}}}^*$
 (Theorem
\ref{open-problem}). The strengthening of G\"{a}rtner-Ellis
theorem (Corollary \ref{GE}) and the solution to the  Ellis
question (Corollary \ref{open-question-Ellis}) are obtained by
taking suitable $\mathcal{S}$.

\begin{proposition}\label{lem-x+}
Let $\lambda_0\in\mathbb{R}$, and assume that $L(\lambda)$  exists
and is finite for all $\lambda$ in an open  interval $G$
containing $\lambda_0$. Then, $(\mu_{\alpha}^{t_\alpha})$ has a
subnet $(\mu_{\gamma}^{t_\gamma})$ such that
 \[l_1^{(\mu_{\gamma}^{t_\gamma})}({L_{\mid G}}'_-(\lambda_0))\le\lambda_0
 {L_{\mid G}}'_-(\lambda_0)-L(\lambda_0)\]
 and
 \[l_1^{(\mu_{\gamma}^{t_\gamma})}({L_{\mid G}}'_+(\lambda_0))\le\lambda_0
 ({L_{\mid G}}'_+(\lambda_0))-L(\lambda_0).\]
 Whence,
 \[l^{(\mu_{\gamma}^{t_\gamma})}_1(x)\le {{L_{\mid G}}}^*(x)\ \ \ \ \ \ \textnormal{for all
 $x\in\textnormal{ran}{L_{\mid G}}_-'\cup\textnormal{ran}{L_{\mid G}}_+'$}.\]
\end{proposition}

\begin{proof}
 Let $G_0$ be an  open interval  such that
$\lambda_0\in G_0\subset\overline{G_0}\subset G $.
 Let $\lambda_1$ and
$\lambda_2$ in $G\verb'\'\{0\}$ such that
$\lambda_1<\lambda<\lambda_2$ for all $\lambda\in G_0$. There
exists $\gamma>1$ such that $\{\gamma \lambda_1,\gamma
\lambda_2\}\subset\mathcal{D}\textnormal{om} (L)$  so that
$h_{\lambda_1}$ and $h_{\lambda_2}$ satisfy (individually) the
tail condition by Lemma 4.3.8 of \cite{dem} (the proof given there
for probability measures works as well for the sub-probability
case). Therefore, for each $\varepsilon>0$ and for each
$i\in\{1,2\}$ there exists $M_{i,\varepsilon}$ such that
\[\limsup\mu_{\alpha}^{t_{\alpha}}(e^{h_{\lambda_i}/t_{\alpha}}1_{\{h_{\lambda_i}
> M_{i,\varepsilon}\}})<\varepsilon.\]
Put $M_\varepsilon=M_{1,\varepsilon}\vee M_{2,\varepsilon}$, and
get
 for each $\lambda\in G_0$,
\[
\int_{\{x:\lambda x> M_\varepsilon\}}e^{\lambda
x/t_\alpha}\mu_\alpha(dx)=\int_{\{x:\lambda x> M_\varepsilon\}\cap
\mathbb{R}_-}e^{\lambda
x/t_\alpha}\mu_\alpha(dx)+\int_{\{x:\lambda x> M_\varepsilon\}\cap
\mathbb{R}_+}e^{\lambda x/t_\alpha}\mu_\alpha(dx)\]
\[\le
\int_{\{x:\lambda_1 x> M_{1,\varepsilon}\}\cap
\mathbb{R}_-}e^{\lambda_1
x/t_\alpha}\mu_\alpha(dx)+\int_{\{x:\lambda_2 x>
M_{2,\varepsilon}\}\cap \mathbb{R}_+}e^{\lambda_2
x/t_\alpha}\mu_\alpha(dx),\] whence
\[\forall \lambda\in G_0,\ \ \ \ \ \ \
 \limsup\mu_{\alpha}^{t_{\alpha}}(e^{h_{\lambda}/t_{\alpha}}1_{\{h_{\lambda}
> M_{\varepsilon}\}})\le\]
\[
\limsup\mu_{\alpha}^{t_{\alpha}}(e^{h_{\lambda_1}/t_{\alpha}}1_{\{h_{\lambda_1}
> M_{1,\varepsilon}\}})
\vee\limsup\mu_{\alpha}^{t_{\alpha}}
(e^{h_{\lambda_2}/t_{\alpha}}1_{\{h_{\lambda_2}
> M_{2,\varepsilon}\}})<\varepsilon.\]
It follows that  $\{h_\lambda:\lambda\in G_0\}$ satisfies the tail
condition for $(\mu_{\alpha}^{t_\alpha})$. Since $L_{\mid G}$ is
continuous and $\overline{G_0}$  compact,  $L_{\mid G_0}$ is
bounded and in particular $\inf_{\lambda\in G_0}L(\lambda)>m$ for
some real $m$. Let $(\mu_{\gamma}^{t_\gamma})$ be a subnet of
$(\mu_{\alpha}^{t_\alpha})$  converging in $\Gamma(X,[0,1])$
(given by Lemma \ref{lem-capa}), put
$\mathcal{T}=\{h_{\lambda}:\lambda\in G_0\}$, and note that all
the hypotheses of
  Theorem \ref{compact-case} hold for $\mathcal{T}$ and
  $(\mu_{\gamma}^{t_\gamma})$, with moreover  $\mathcal{T}\subset\mathcal{C}_{\mathcal{K}}(X)$.
If $\lambda_0\neq 0$ (say $\lambda_0>0$), then $\lambda_1$ and
$\lambda_2$ can be chosen such that
$0<\lambda_1<\lambda<\lambda_2$ for all $\lambda\in G_0$.
 Since for each real $M\ge m$, there is a compact $K_M$ such that
$\bigcup_{\lambda\in G_0}\{m\le h_\lambda\le M\}\subset K_M$, by
Theorem \ref{compact-case} (a)
 we get  a compact $K$ such that
 \begin{equation}\label{lem-x+-eq1}
 L(\lambda)=\sup_{x\in K}
  \{\lambda x-l^{(\mu_{\gamma}^{t_\gamma})}_1(x)\}\ \ \ \ \ \
  \textnormal{for all $\lambda\in G_0$}.
  \end{equation}
  If $\lambda_0=0$, then  $(\mu_\alpha)$ (resp. $(\mu_\gamma)$) is
   exponentially tight with respect to
 $(t_{\alpha})$ (resp. $(t_{\gamma})$), and we apply
 Theorem \ref{compact-case} (b)  to get (\ref{lem-x+-eq1}).
  Therefore,
  for each $\lambda\in G_0$ there exists $x_\lambda\in K$ such that
  $L(\lambda)=
  \lambda x_\lambda-l^{(\mu_{\gamma}^{t_\gamma})}_1(x_\lambda)$.
  Put $x={L_{\mid G}}_+'(\lambda_0)$, and
let $(x_{\lambda'+\lambda_0})$ be a subnet of
$(x_{\lambda+\lambda_0})_{\lambda+\lambda_0\in G_0,\lambda>0}$.
 Since  $x_{\lambda+\lambda_0}\in K$ for all
$\lambda+\lambda_0\in G_0$, $(x_{\lambda'+\lambda_0})$  has a
subnet $(x_{\lambda''+\lambda_0})$ converging to some point
$x''\in K$  when $\lambda''\rightarrow 0^+$, so that
\[x=\lim_{\lambda''\rightarrow
0^+}\frac{L(\lambda''+\lambda_0)-L(\lambda_0)}{\lambda''}=
\lim_{\lambda''\rightarrow 0^+}\frac{
(\lambda''+\lambda_0)x_{\lambda''+\lambda_0}
-l^{(\mu_{\gamma}^{t_\gamma})}_1(x_{\lambda''+\lambda_0})-L(\lambda_0)}{\lambda''}\]
\[= x''+\lim_{\lambda''\rightarrow 0^+}\frac{\lambda_0
x_{\lambda''+\lambda_0}-l^{(\mu_{\gamma}^{t_\gamma})}_1(x_{\lambda''+\lambda_0})-L(\lambda_0)}{\lambda''},\]
which implies $x''=x$ and
\[0=\lim_{\lambda''\rightarrow 0^+}\lambda_0
x_{\lambda''+\lambda_0}-l^{(\mu_{\gamma}^{t_\gamma})}_1(x_{\lambda''+\lambda_0})-L(\lambda_0)\le\lambda_0
x-l^{(\mu_{\gamma}^{t_\gamma})}_1(x)-L(\lambda_0),
\]which proves the assertion concerning ${L_{\mid G}}_+'(\lambda_0)$. A similar proof
works for ${L_{\mid G}}_-'(\lambda_0)$.
\end{proof}

\begin{theorem}\label{open-problem}
 Let $\mathcal{S}\subset\mathcal{C}(X)$
and $G\subset X$ be a nonempty open interval  such
 that
 $\mathcal{S}\supset
\{h_\lambda:\lambda\in G\}$, and assume that
 $\Lambda(h)$ exists for all $h\in\mathcal{S}$ with
$L(\lambda)$  finite for all $\lambda\in G$.
\begin{itemize}
\item[(a)]
If \begin{equation}\label{open-problem-eq0.01}
 \textnormal{ran}{{L_{\mid G}}_-'}\cup\textnormal{ran}{{L_{\mid G}}_+'}
 \supset\mathcal{D}\textnormal{om}(l_0)\cap\{l_1>-\overline{\Lambda}(0)\},
\end{equation}
then $(\mu_\alpha)$ satisfies a vague large deviation principle
with powers $(t_{\alpha})$ and rate function  $J$ satisfying
 \begin{equation}\label{open-problem-eq0.02}
J(x)={L_{\mid G}}^*(x)={\Lambda_{\mid \mathcal{S}}}^*(x)\ \ \ \ \
\ \textnormal{for all
$x\in\mathcal{D}\textnormal{om}(J)\cap\{J>-\overline{\Lambda}(0)\}$}.
\end{equation}
If moreover $0\in G$, then the principle is narrow and
\begin{equation}\label{open-problem-eq0.03}
J(x)={L_{\mid G}}^*(x)={\Lambda_{\mid \mathcal{S}}}^*(x)\ \ \ \ \
\ \textnormal{for all $x\in\mathcal{D}\textnormal{om}(J)$}.
\end{equation}
\item[(b)]
If \begin{equation}\label{open-problem-eq0.04}
 \textnormal{ran}{{L_{\mid G}}_-'}\cup\textnormal{ran}{{L_{\mid G}}_+'}
 \supset\mathcal{D}\textnormal{om}(l_0),
\end{equation}
then $(\mu_\alpha)$ satisfies a vague large deviation principle
with powers $(t_{\alpha})$ and rate function  $J$ satisfying
 \begin{equation}\label{open-problem-eq0.05}
J(x)={L_{\mid G}}^*(x)={\Lambda_{\mid \mathcal{S}}}^*(x)\ \ \ \ \
\ \textnormal{for all $x\in\mathcal{D}\textnormal{om}(J)$}.
\end{equation}
If moreover $0\in G$, then the principle is narrow.
\item[(c)]
If
\begin{equation}\label{open-problem-eq0.06}
 \textnormal{ran}{{L_{\mid G}}_-'}\cup\textnormal{ran}{{L_{\mid G}}_+'}
 \supset\mathcal{D}\textnormal{om}({\Lambda_{\mid
\mathcal{S}}}^*)\cap\{l_1>-\overline{\Lambda}(0)\},
\end{equation} then $(\mu_\alpha)$ satisfies a vague large deviation principle with
powers $(t_{\alpha})$ and rate function $J$  satisfying
\begin{equation}\label{open-problem-eq0.07}
J(x)={\Lambda_{\mid \mathcal{S}}}^*(x)\ \ \ \ \ \ \textnormal{for
all $x\in\{J>-\overline{\Lambda}(0)\}$},
\end{equation}
and
\begin{equation}\label{open-problem-eq0.08}
 J(x)={L_{\mid G}}^*(x)\ \ \ \ \ \ \textnormal{for
all $x\in\mathcal{D}\textnormal{om}({\Lambda_{\mid
\mathcal{S}}}^*)\cap\{J>-\overline{\Lambda}(0)\}$}.
\end{equation}
 If moreover $0\in G$, then the principle is narrow with $J={\Lambda_{\mid
 \mathcal{S}}}^*$ satisfying
\begin{equation}\label{open-problem-eq0.09}
J(x)={L_{\mid G}}^*(x)\ \ \ \ \ \ \textnormal{for all
$x\in\mathcal{D}\textnormal{om}(J)$}.
\end{equation}
\item[(d)]
 If \begin{equation}\label{open-problem-eq0.1}
 \textnormal{ran}{{L_{\mid G}}_-'}\cup\textnormal{ran}{{L_{\mid G}}_+'}
 \supset\mathcal{D}\textnormal{om}({\Lambda_{\mid
\mathcal{S}}}^*),
\end{equation}
then $(\mu_\alpha)$ satisfies a vague large deviation principle
with powers $(t_{\alpha})$ and rate function $J={\Lambda_{\mid
 \mathcal{S}}}^*$ satisfying
\begin{equation}\label{open-problem-eq0.2}
 J(x)={L_{\mid G}}^*(x)\ \ \ \ \ \ \textnormal{for
all $x\in\mathcal{D}\textnormal{om}(J)$}.
\end{equation}
If moreover $0\in G$, then the principle is narrow.
\item[(e)] If  $l_0$ is  proper convex, then $(a)$ (resp. $(b)$)  holds verbatim replacing the symbol
$\mathcal{D}\textnormal{om}$ by
$\textnormal{int}\mathcal{D}\textnormal{om}$ in
(\ref{open-problem-eq0.01}), (\ref{open-problem-eq0.02}),
(\ref{open-problem-eq0.03}) (resp. (\ref{open-problem-eq0.04}),
(\ref{open-problem-eq0.05})).
\item[(f)] If  ${\Lambda_{\mid \mathcal{S}}}^*$ is  proper convex
and lower semi-continuous, then $(c)$ (resp. $(d)$)  holds
verbatim replacing the symbol $\mathcal{D}\textnormal{om}$ by
$\textnormal{int}\mathcal{D}\textnormal{om}$ in
(\ref{open-problem-eq0.06}),  (\ref{open-problem-eq0.08}),
(\ref{open-problem-eq0.09})) (resp. (\ref{open-problem-eq0.1}),
(\ref{open-problem-eq0.2})).
\end{itemize}
\end{theorem}

\begin{proof}
 For all $h\in\mathcal{S}$ and all $x\in X$ we have by Theorem
 \ref{lim} (since $\Lambda(h)\ge\overline{\Lambda}(h1_{\{h\le M\}}+(-\infty)1_{\{h> M\}})$ for all reals $M$),
\[\Lambda(h)-h(x)\ge\sup_{M\in\mathbb{R}}\sup_{\{h\le M\}}\{h(y)-l_0(y)\}-h(x)\ge\sup_{y\in X}\{h(y)-l_0(y)\}-h(x)\ge
-l_0(x),\]so that
\begin{equation}\label{open-problem-eq1}
{L_{\mid G}}^*(x)\le{\Lambda_{\mid \mathcal{S}}}^*(x)\le l_0(x)\ \
\ \ \ \textnormal{for all $x\in X$}.
\end{equation}
Assume that (\ref{open-problem-eq0.01}) holds, and  let
$(\mu_{\beta}^{t_\beta})$ be a subnet of
$(\mu_{\alpha}^{t_\alpha})$. By Proposition \ref{lem-x+} applied
to $(\mu_{\beta}^{t_\beta})$ in place of
$(\mu_{\alpha}^{t_\alpha})$, $(\mu_{\beta}^{t_\beta})$ has a
subnet $(\mu_{\gamma}^{t_\gamma})$ such that
\begin{equation}\label{open-problem-eq1.0}
l^{(\mu_{\gamma}^{t_\gamma})}_1(x) \le{{L_{\mid G}}}^*(x)\ \ \ \ \
\ \textnormal{for all
 $x\in\mathcal{D}\textnormal{om}(l_0)\cap\{l_1>-\overline{\Lambda}(0)\}$}.\end{equation}
 Since
 \begin{equation}\label{open-problem-eq1.1}
 l_0\le l^{(\mu_{\gamma}^{t_\gamma})}_0\le l^{(\mu_{\gamma}^{t_\gamma})}_1,
 \end{equation}
 (\ref{open-problem-eq1}) and  (\ref{open-problem-eq1.0}) imply
\begin{equation}\label{open-problem-eq2}
l^{(\mu_{\gamma}^{t_\gamma})}_0(x)=l^{(\mu_{\gamma}^{t_\gamma})}_1(x)={L_{\mid
G}}^*(x)={\Lambda_{\mid \mathcal{S}}}^*(x)=l_0(x)\ \ \ \ \ \
\textnormal{for all
$x\in\mathcal{D}\textnormal{om}(l_0)\cap\{l_1>-\overline{\Lambda}(0)\}$}.
\end{equation} If
$x\not\in\mathcal{D}\textnormal{om}(l_0)$, then
$l^{(\mu_{\gamma}^{t_\gamma})}_0(x)=l^{(\mu_{\gamma}^{t_\gamma})}_1(x)=+\infty$
by   (\ref{open-problem-eq1.1}). If $l_1(x)\le
-\overline{\Lambda}(0)$, then
\[l^{(\mu_{\gamma}^{t_\gamma})}_0(x)=l^{(\mu_{\gamma}^{t_\gamma})}_1(x)=l_0(x)=l_1(x)=-\overline{\Lambda}(0).\]
Therefore,
 $l^{(\mu_{\gamma}^{t_\gamma})}_0(x)=l^{(\mu_{\gamma}^{t_\gamma})}_1(x)$ for all $x\in X$.
By Lemma \ref{lem-prelim} applied to $(\mu_{\gamma}^{t_\gamma})$,
$(\mu_{\gamma})$ satisfies a vague large deviation principle with
powers $(t_\gamma)$ and rate function
\begin{equation}\label{open-problem-eq3}
J(x)=\left\{
\begin{array}{ll}
{\Lambda_{\mid \mathcal{S}}}^* & \textnormal{if} \
x\in\mathcal{D}\textnormal{om}(l_0)\cap\{l_1>-\overline{\Lambda}(0)\}\\
\\
-\overline{\Lambda}(0) & \textnormal{if}\  l_1(x)\le -\overline{\Lambda}(0)\\
\\
+\infty & \textnormal{if}\
x\not\in\mathcal{D}\textnormal{om}(l_0).
 \end{array}
 \right.
\end{equation}
By Lemma \ref{lem-capa} (b), $(\mu_{\gamma}^{t_\gamma})$ converges
to $c_J$ in $\Gamma(X,[0,1])$. Since $(\mu_{\beta}^{t_\beta})$ is
arbitrary, we have proved that any subnet of
$(\mu_{\alpha}^{t_\alpha})$ has a subnet converging vaguely to
$c_J$.  By Lemma \ref{lem-capa} (a), it follows  that
$(\mu_{\alpha}^{t_\alpha})$ converges vaguely to $c_J$, which
proves  the first assertion of (a) ((\ref{open-problem-eq0.02})
follows from (\ref{open-problem-eq2}) and
(\ref{open-problem-eq3}),  since $J=l_0=l_1$). If $0\in G$, then
(\ref{open-problem-eq0.03}) follows from (\ref{open-problem-eq1})
and (\ref{open-problem-eq3}) since $-L(0)\le{L_{\mid G}}^*$, and
the principle is narrow by exponential tightness. The proofs of
(b),(c),(d) are similar. Assume that $l_0$ is proper convex, and
\[
 \textnormal{ran}{{L_{\mid G}}_-'}\cup\textnormal{ran}{{L_{\mid G}}_+'}
 \supset\textnormal{int}\mathcal{D}\textnormal{om}(l_0)\cap\{l_1>-\overline{\Lambda}(0)\}.
\]
 In the same way as above we get
 \begin{equation}\label{open-problem-eq4}
l^{(\mu_{\gamma}^{t_\gamma})}_0(x)=l^{(\mu_{\gamma}^{t_\gamma})}_1(x)={L_{\mid
G}}^*(x)={\Lambda_{\mid \mathcal{S}}}^*(x)=l_0(x)\ \ \ \ \ \
\textnormal{for all
$x\in\textnormal{int}\mathcal{D}\textnormal{om}(l_0)\cap\{l_1>-\overline{\Lambda}(0)\}$}.
\end{equation}
Suppose that $l^{(\mu_{\gamma}^{t_\gamma})}_1(x)>l_0(x)$ for some
$x\in\{l_1>-\overline{\Lambda}(0)\}$. Since $l_1$ and
$l^{(\mu_{\gamma}^{t_\gamma})}_1$ are lower semi-continuous, there
is an open set $G_0$ containing $x$ such that
\[\inf_{G_0\cap\{l_1>-\overline{\Lambda}(0)\}}
l^{(\mu_{\gamma}^{t_\gamma})}_1>\inf_{G_0\cap\{l_1>-\overline{\Lambda}(0)\}}
 l_0=\inf_{G_0\cap\{l_1>-\overline{\Lambda}(0)\}\cap\textnormal{int}\mathcal{D}\textnormal{om}(
l_0)} l_0,\]
 where the
equality follows from Lemma \ref{conv}  applied to $l_0$ and
$G_0\cap\{l_1>-\overline{\Lambda}(0)\}$. Then, there exists $y\in
G_0\cap\{l_1>-\overline{\Lambda}(0)\}\cap\textnormal{int}\mathcal{D}\textnormal{om}(l_0)$
such that $l^{(\mu_{\gamma}^{t_\gamma})}_1(y)>l_0(y)$, which
contradicts (\ref{open-problem-eq4}). We then have
$l^{(\mu_{\gamma}^{t_\gamma})}_1(x)\le l_0(x)$ for all
$x\in\{l_1>-\overline{\Lambda}(0)\}$, and by
(\ref{open-problem-eq1.1}),
\[
l^{(\mu_{\gamma}^{t_\gamma})}_0(x)=l^{(\mu_{\gamma}^{t_\gamma})}_1(x)=l_0(x)\
\ \ \ \ \ \textnormal{for all
$x\in\{l_1>-\overline{\Lambda}(0)\}$}.
\]
Since \[
l^{(\mu_{\gamma}^{t_\gamma})}_0(x)=l^{(\mu_{\gamma}^{t_\gamma})}_1(x)=
l_0(x)=l_1(x)=-\overline{\Lambda}(0) \ \ \ \ \ \ \textnormal{for
all $x\in\{l_1\le-\overline{\Lambda}(0)\}$},
\]
 it follows as above that $(\mu_{\alpha}^{t_\alpha})$ converges
vaguely to $c_J$, with $J=l_0=l_1$ satisfying by
(\ref{open-problem-eq4}),
\begin{equation}\label{open-problem-eq6}
J(x)={L_{\mid G}}^*(x)={\Lambda_{\mid \mathcal{S}}}^*(x)\ \ \ \ \
\ \textnormal{for all
$x\in\textnormal{int}\mathcal{D}\textnormal{om}(J)\cap\{J>-\overline{\Lambda}(0)\}$}.
\end{equation}
If $0\in G$, then $-L(0)\le{L_{\mid G}}^*$,  and by
(\ref{open-problem-eq1})  and (\ref{open-problem-eq6})  we get
\[
J(x)={L_{\mid G}}^*(x)={\Lambda_{\mid \mathcal{S}}}^*(x)\ \ \ \ \
\ \textnormal{for all
$x\in\textnormal{int}\mathcal{D}\textnormal{om}(J)$}.
\]
This proves the assertion of (e) concerning (a); the one
concerning (b) is proved similarly.
 Assume that ${\Lambda_{\mid
\mathcal{S}}}^*$ is proper convex lower semi-continuous,  and
\[
 \textnormal{ran}{{L_{\mid G}}_-'}\cup\textnormal{ran}{{L_{\mid G}}_+'}
 \supset\textnormal{int}\mathcal{D}\textnormal{om}({\Lambda_{\mid
\mathcal{S}}}^*)\cap\{l_1>-\overline{\Lambda}(0)\}.
\]
 As above we get
 \begin{equation}\label{open-problem-eq8}
l^{(\mu_{\gamma}^{t_\gamma})}_0(x)=l^{(\mu_{\gamma}^{t_\gamma})}_1(x)={L_{\mid
G}}^*(x)={\Lambda_{\mid \mathcal{S}}}^*(x)\ \ \ \ \ \
\textnormal{for all
$x\in\textnormal{int}\mathcal{D}\textnormal{om}({\Lambda_{\mid
\mathcal{S}}}^*)\cap\{l_1>-\overline{\Lambda}(0)\}$}.
\end{equation}
The same reasoning as in the proof of (e) (with ${\Lambda_{\mid
\mathcal{S}}}^*$ in place of $l_0$) gives
$l^{(\mu_{\gamma}^{t_\gamma})}_1(x)\le{\Lambda_{\mid
\mathcal{S}}}^*(x)$ for all $x\in\{l_1>-\overline{\Lambda}(0)\}$,
and by (\ref{open-problem-eq1}),
\begin{equation}\label{open-problem-eq9}
l^{(\mu_{\gamma}^{t_\gamma})}_0(x)=l^{(\mu_{\gamma}^{t_\gamma})}_1(x)={\Lambda_{\mid
\mathcal{S}}}^*(x)=l_0(x)\ \ \ \ \ \ \textnormal{for all
$x\in\{l_1>-\overline{\Lambda}(0)\}$}.
\end{equation}
Since \[
l^{(\mu_{\gamma}^{t_\gamma})}_0(x)=l^{(\mu_{\gamma}^{t_\gamma})}_1(x)=-\overline{\Lambda}(0)
\ \ \ \ \ \ \textnormal{for all
$x\in\{l_1\le-\overline{\Lambda}(0)\}$},
\]
 it follows as above that $(\mu_{\alpha}^{t_\alpha})$ converges
vaguely to $c_J$, with $J$ satisfying (\ref{open-problem-eq0.07}).
Since $J=l_1$, (\ref{open-problem-eq8}) gives
 \begin{equation}\label{open-problem-eq10}
J(x)={L_{\mid G}}^*(x)\ \ \ \ \ \ \textnormal{for all
$x\in\textnormal{int}\mathcal{D}\textnormal{om}({\Lambda_{\mid
\mathcal{S}}}^*)\cap\{J>-\overline{\Lambda}(0)\}$}.
\end{equation}
Since $0\in G$ implies $-L(0)\le{L_{\mid G}}^*$, by
(\ref{open-problem-eq1}), (\ref{open-problem-eq9}),
(\ref{open-problem-eq10}), we obtain $J={\Lambda_{\mid
\mathcal{S}}}^*$ and
\[
J(x)={L_{\mid G}}^*(x)\ \ \ \ \ \ \textnormal{for all
$x\in\textnormal{int}\mathcal{D}\textnormal{om}({\Lambda_{\mid
\mathcal{S}}}^*)$}.
\]
This proves the assertion of (f) concerning (c); the one
concerning (d) is proved similarly.
\end{proof}

 The standard
G\"{a}rtner-Ellis theorem deals with the case where $(\mu_\alpha)$
is a sequence of Borel probability measures; it states that if
$L(\lambda)$ exists  for all reals $\lambda$,
 $L$ is lower
 semi-continuous
  essentially  smooth and
  $0\in\textnormal{int}\mathcal{D}\textnormal{om}(L)$,
 then $(\mu_\alpha)$ satisfies a  large deviation principle with
powers $(t_{\alpha})$ and rate function $L^*$ (\cite{dem}, Theorem
2.3.6, \cite{gar}, \cite{ell3}). A stronger version has been given
by O' Brien (\cite{ob}, Theorem 5.1): if $L(\lambda)$ exists and
is finite for all $\lambda$  in a nonempty open interval $G$ and
if $\widehat{L_{\mid G}}$ is essentially smooth,
 then $(\mu_\alpha)$ satisfies a vague large deviation principle with
powers $(t_{\alpha})$ and rate function ${L_{\mid G}}^*$; if
moreover $0\in G$, then the principle is narrow. The former
version is recovered by taking
$G=\textnormal{int}\mathcal{D}\textnormal{om}(L)$ (the hypotheses
implying $L^*={L_{\mid G}}^*$ with $\widehat{L_{\mid G}}$
essentially smooth). The improvements consists in the obtention of
the vague large deviations, and in the fact that $L$ in not
assumed to exist out $G$ (even when  $L$ exists on $X$, it is not
assumed to be lower semi-continuous).

The following corollary  summarizes the case where
$\mathcal{S}=\{h_\lambda:\lambda\in G\}$ in Theorem
\ref{open-problem}, and where large deviations hold with rate
function ${L_{\mid G}}^*$ ($={\Lambda_{\mid \mathcal{S}}}^*$). It
strengthens the O' Brien's version of G\"{a}rtner-Ellis theorem by
obtaining the same conclusions, with the essential smoothness
hypothesis replaced by the weaker condition (\ref{GE-eq0}) (or
(\ref{GE-eq0.1}) when $0\in G$);
 in particular, there is no
differentiability assumption. Furthermore, it works for general
nets of Radon sub-probability measures.

\begin{corollary}\label{GE}
We assume that $L(\lambda)$ exists and is finite for all $\lambda$
 in a nonempty open interval $G\subset X$.
\begin{itemize}
\item[(a)] If
\begin{equation}\label{GE-eq0}
\textnormal{ran}{L_{\mid G}}'_-\cup\textnormal{ran}{L_{\mid
G}}'_+\supset\textnormal{int}\mathcal{D}\textnormal{om}({L_{\mid
G}}^*),
\end{equation}then
 $(\mu_\alpha)$ satisfies a vague large deviation principle with
powers $(t_{\alpha})$ and rate function ${L_{\mid G}}^*$. The
condition (\ref{GE-eq0}) is satisfied in particular when
$\widehat{L_{\mid G}}$ is essentially smooth.
\item[(b)]  If $0\in G$ and
\begin{equation}\label{GE-eq0.1}
\textnormal{ran}{L_{\mid G}}'_-\cup\textnormal{ran}{L_{\mid
G}}'_+\supset\textnormal{int}\mathcal{D}\textnormal{om}({L_{\mid
G}}^*)\cap\{l_1>-L(0)\},
\end{equation}then
 $(\mu_\alpha)$ satisfies a large deviation principle with
powers $(t_{\alpha})$ and rate function ${L_{\mid G}}^*$.
\end{itemize}
\end{corollary}

\begin{proof}
(b) and the first assertion of (a) follow from Theorem
\ref{open-problem} (f) with $\mathcal{S}=\{h_\lambda:\lambda\in
G\}$. Assume that $\widehat{L_{\mid G}}$  is essentially smooth.
Extend $L_{\mid G}$ by continuity to a convex function $L_{\mid
\overline{G}}$ on $\overline{G}$, so that $\widehat{L_{\mid
\overline{G}}}$ is  a proper  convex lower semi-continuous
function on $X$ with
$G=\textnormal{int}\mathcal{D}\textnormal{om}(\widehat{L_{\mid
\overline{G}}})$; moreover, $\widehat{L_{\mid \overline{G}}}$ is
essentially smooth. By   Theorem 26.1 and Corollary 26.4.1 of
\cite{roc},
\begin{equation}\label{GE-eq1}
\textnormal{ran}{L_{\mid
\overline{G}}}'\supset\textnormal{int}\mathcal{D}\textnormal{om}({L_{\mid
\overline{G}}}^*),
\end{equation}
which gives (\ref{GE-eq0})  since $\textnormal{ran}{L_{\mid
\overline{G}}}'=\textnormal{ran}{L_{\mid {G}}}'$  and ${L_{\mid
\overline{G}}}^*=L_{\mid G}^*$.
\end{proof}

The solution to the Ellis question (with in fact weaker
hypotheses) is a direct consequence of Theorem \ref{open-problem}
$(c)$,  by taking
$\mathcal{S}=\{h_{\lambda,\nu}:(\lambda,\nu)\in\mathbb{R}^2\}$.

\begin{corollary}\label{open-question-Ellis}
Put
$\mathcal{S}=\{h_{\lambda,\nu}:(\lambda,\nu)\in\mathbb{R}^2\}$,
and assume that $\Lambda(h_{\lambda,\nu})$ exists for all
$(\lambda,\nu)\in\mathbb{R}^2$ and is finite for all pairs
$(\lambda,\lambda)$ with $\lambda$ in some open interval $G$
containing $0$. If
  $\textnormal{ran}{L_{\mid
G}}'_-\cup\textnormal{ran}{L_{\mid
G}}'_+\supset\mathcal{D}\textnormal{om}({\Lambda_{\mid
\mathcal{S}}}^*)\cap\{l_1>-L(0)\}$, then $(\mu_\alpha)$ satisfies
a  large deviation principle with powers $(t_{\alpha})$ and rate
function $J={\Lambda_{\mid \mathcal{S}}}^*$. Moreover,
\[
J(x)={L_{\mid G}}^*(x)\ \ \ \ \ \ \textnormal{for all
$x\in\mathcal{D}\textnormal{om}(J)$}.
\]
\end{corollary}

The following example is often cited as a typical case not covered
by the G\"{a}rtner-Ellis theorem (\cite{ell1}, \cite{ob}).

\begin{example}\label{GE-ex}
Consider the sequence  $(\mu_n^{1/n})$ where
$\mu_n\{-1\}=\mu_n\{1\}=\frac{1}{2}$ for all $n\in\mathbb{N}$.
Then $L(\lambda)=|\lambda|$ for all reals $\lambda$. Take
$\mathcal{S}=\{h_{\lambda,\nu}:(\lambda,\nu)\in\mathbb{R}^2\}$ and
compute
\[\Lambda(h_{\lambda,\nu})=-\lambda\vee \nu\ \ \ \ \ \ \ \ \textnormal{for all $(\lambda,\nu)\in\mathbb{R}^2$},\]
whence
\begin{displaymath}
{\Lambda_{\mid \mathcal{S}}}^*(x)=\left\{
\begin{array}{ll}
0 & \textnormal{if}\  |x|=1
\\
+\infty & \textnormal{if}\  |x|\neq 1.
\end{array}
\right.
\end{displaymath} Then,
 $\textnormal{ran}{L_-'}\cup\textnormal{ran}{L_+'}=
 \{-1,1\}\supset\mathcal{D}\textnormal{om}({\Lambda_{\mid
\mathcal{S}}}^*)$, and
 by Corollary \ref{open-question-Ellis},  $(\mu_n)$
satisfies a large deviation principle with powers $(1/n)$ and
 rate function $J={\Lambda_{\mid \mathcal{S}}}^*$. Since
\begin{displaymath}
L^*(x)=\left\{
\begin{array}{ll}
0 & \textnormal{if}\  |x|\le 1
\\
+\infty & \textnormal{if}\  |x|> 1,
\end{array}
\right.
\end{displaymath}
we have  $J(x)=L^*(x)$ for all
$x\in\{-1,1\}=\mathcal{D}\textnormal{om}(J)$. Note that for any
nonempty open set $G\subset]-1,1[$,
\[\textnormal{ran}{L_{\mid G}}'_-\cup\textnormal{ran}{L_{\mid G}}'_+\not\supset
\textnormal{int}\mathcal{D}\textnormal{om}(L_{\mid
G}^*)\cap\{J>0\}\supset]-1,1[,\] and the condition
(\ref{GE-eq0.1}) of Corollary \ref{GE} does not hold.
\end{example}

The following example exhibits a situation with convex rate
function, where both above corollaries do not  work; we then apply
 theorem \ref{open-problem} with another set $\mathcal{S}$.

\begin{example}\label{ex-Dem-Zei}
Consider the net
$(\mu_\varepsilon^{\varepsilon})_{\varepsilon>0}$, where
$\mu_\varepsilon$ is the  probability measure on $X$ defined by
$\mu_\varepsilon(0)=1-2p_\varepsilon$,
$\mu_\varepsilon(-\varepsilon\log
p_\varepsilon)=\mu_\varepsilon(\varepsilon\log
p_\varepsilon)=p_\varepsilon$, and assume that $\lim
\varepsilon\log p_\varepsilon=-\infty$. Put
$Q_n(x)=n|x|e^{-|x|}-x$ for all $n\in\mathbb{N}$ and all $x\in X$,
and take
$\mathcal{S}=\{Q_n:n\in\mathbb{N}\}\cup\{h_\lambda:\lambda\in]-1,1[\}$.
Easy calculations give $\Lambda(Q_n)=0$ for all $n\in\mathbb{N}$,
and \begin{displaymath}
 L(\lambda)=\left\{
\begin{array}{ll}
0 & \ \ \ \ \ \ \textnormal{if $|\lambda|\le 1$}
\\
+\infty & \ \ \ \ \ \ \textnormal{if $|\lambda|>1$},
\end{array}
\right.
\end{displaymath}
so that
\[{L_{\mid ]-1,1[}}^*(x)=L^*(x)=|x|\ \ \ \ \ \ \textnormal{for all
$x\in X$},\] and
\begin{displaymath}
{\Lambda_{\mid
\mathcal{S}}}^*(x)=\sup_{n\in\mathbb{N}}\{Q_n(x)-\Lambda(Q_n)\}\vee
{L_{\mid ]-1,1[}}^*(x)=\left\{
\begin{array}{ll}
0 & \textnormal{if $x=0$}
\\
+\infty & \textnormal{otherwise}.
\end{array}
\right.
\end{displaymath}
Then, $\textnormal{ran}{L_{\mid
]-1,1[}}'=\{0\}\supset\mathcal{D}\textnormal{om}({\Lambda_{\mid
\mathcal{S}}}^*)$, and by Theorem  \ref{open-problem} $(d)$,
$(\mu_\varepsilon)$ satisfies a large deviation principle with
powers $(\varepsilon)_{\varepsilon>0}$ and rate function
$J={\Lambda_{\mid \mathcal{S}}}^*$. Note that $J$ is convex but
$J\neq L^*$ (however, $J$ coincides with $L^*$ on
$\mathcal{D}\textnormal{om}(J)$); in particular, $L$ is not
essentially smooth and the G\"{a}rtner-Ellis theorem does not
work. Furthermore,  for any nonempty open set $G\subset]-1,1[$,
\[\{0\}=\textnormal{ran}{L_{\mid G}}'\not\supset
\textnormal{int}\mathcal{D}\textnormal{om}({L_{\mid G
}}^*)\cap\{J>0\}\supset X\verb'\'\{0\}\] and the condition
(\ref{GE-eq0.1}) of Corollary  \ref{GE} does not hold either. We
observe also that Corollary \ref{open-question-Ellis} does not
apply; indeed,   the set
$\{h_{\lambda,\nu}:(\lambda,\nu)\in\mathbb{R}^2\}$ is not suitable
 since
\begin{displaymath}
\Lambda(h_{\lambda,\nu})=\left\{
\begin{array}{ll}
0 & \textnormal{if $\lambda\ge -1$ and $\nu\le 1$}
\\
+\infty & \textnormal{otherwise}
\end{array}
\right.
\end{displaymath}
 gives
 ${\Lambda_{\mid \{h_{\lambda,\nu}:(\lambda,\nu)\in\mathbb{R}^2\}}}^*(x)=L^*(x)$ for all
 $x\in X$.
\end{example}

\end{document}